\documentstyle[11pt]{article}
\title{Some examples of rigid representations\footnote{Research partially 
supported by INTAS grant 97-1644}}
\author{Vladimir Petrov Kostov\\ \\ \hspace{7cm}{\sl To my mother}} 
\date{}
\newtheorem{tm}{Theorem}
\newtheorem{lm}[tm]{Lemma}
\newtheorem{cor}[tm]{Corollary}
\newtheorem{prop}[tm]{Proposition}
\newtheorem{rem}[tm]{Remark}
\newtheorem{ex}[tm]{Example}
\newtheorem{defi}[tm]{Definition}
\begin{document}
\maketitle

\begin{abstract}
Consider the Deligne-Simpson problem: {\em give necessary and 
sufficient conditions for the choice of the conjugacy classes 
$C_j\subset GL(n,{\bf C})$ (resp. $c_j\subset gl(n,{\bf C})$) so that there 
exist irreducible $(p+1)$-tuples of matrices $M_j\in C_j$ (resp. $A_j\in c_j$) 
satisfying the equality $M_1\ldots M_{p+1}=I$ 
(resp. $A_1+\ldots +A_{p+1}=0$)}.   
The matrices $M_j$ and $A_j$ are interpreted as monodromy operators  
and as matrices-residua of fuchsian systems on Riemann's 
sphere. 

We give new examples of existence of such $(p+1)$-tuples of matrices $M_j$ 
(resp. $A_j$)
which are {\em rigid}, i.e. unique up to conjugacy once the classes $C_j$ 
(resp. $c_j$) are fixed. For rigid representations the sum of the dimensions 
of the classes $C_j$ (resp. $c_j$) equals $2n^2-2$.  
\end{abstract}

\section{Fuchsian linear systems and the Deligne-Simpson problem}

Consider the {\em fuchsian} system (i.e. with logarithmic poles) of 
$n$ linear differential equations 

\[ {\rm d}X/{\rm d}t=(\sum _{j=1}^{p+1}A_j/(t-a_j))X~,~t\in {\bf C}P^1=
{\bf C}\cup \infty \] 
$A_j\in gl(n,{\bf C})$ being its {\em matrices-residua}. Assume that it has no 
pole at infinity, i.e.

\begin{equation}\label{A_j}
A_1+\ldots +A_{p+1}=0~.
\end{equation}

Fix a base point $a_0\in S:={\bf C}P^1\backslash \{ a_1,\ldots ,
a_{p+1}\}$ 
and the value $B\in GL(n,{\bf C})$ of the solution $X$ for $t=a_0$. 
Every pole $a_j$ defines a conjugacy class $\gamma _j$ in the fundamental 
group $\pi _1(S,a_0)$.
The class $\gamma _j$ is represented 
by a closed contour consisting of a segment $[a_0,a_j']$ 
(where the point $a_j'$ is close to $a_j$), 
of a circumference centered at $a_j$ and of radius $|a_j-a_j'|$
(containing inside no pole of the system other 
than $a_j$ and circumventing $a_j$ counterclockwise) and of the segment 
$[a_j',a_0]$. One enumerates the segments so that the index 
increases when one turns around $a_0$ clockwise. 

Hence, 
$\pi _1(S,a_0)$ admits the presentation 

\[ \left< \gamma _1,\ldots ,\gamma_ {p+1}|\gamma _{p+1}
\ldots \gamma _1=e\right>~.\] 
The {\em monodromy operator} $M_j$ 
defined by the class $\gamma _j$ is the one mapping 
the solution with initial data $X|_{t=a_0}=B$ onto the value at $a_0$ of 
its analytic continuation 
along the contour defining $\gamma _j$ (i.e. $X\mapsto XM_j$). The monodromy 
operators of the system generate its {\em monodromy group} which is an  
antirepresentation $\pi _1(S,a_0)\rightarrow GL(n,{\bf C})$ 
because the monodromy operator corresponding to the 
class $\gamma _i\gamma _j$ equals $M_jM_i$. Thus for the matrices $M_j$ 
one has 

\begin{equation}\label{M_j}
M_1\ldots M_{p+1}=I
\end{equation}

\begin{rem}\label{2pii}
{\em If there are no non-zero integer differences 
between the eigenvalues of $A_j$, then the operator $M_j$ is conjugate to 
$\exp (2\pi iA_j)$.} 
\end{rem}

\begin{rem}\label{regular}
{\em Fuchsian systems are a particular case of {\em regular} systems, i.e. 
linear systems whose solutions when restricted to sectors centered at the 
poles $a_j$ grow no faster than some power of the distance to the pole $a_j$. 
Their monodromy groups are defined in the same way.} 
\end{rem}

The {\em Deligne-Simpson problem (DSP)} is formulated like this: {\em give 
necessary and sufficient conditions for the choice of the conjugacy classes 
$C_j\subset GL(n,{\bf C})$ (resp. $c_j\subset gl(n,{\bf C})$) so that 
there exist irreducible $(p+1)$-tuples of matrices $M_j\in C_j$ satisfying 
(\ref{M_j}) (resp. of matrices $A_j\in c_j$ satisfying 
(\ref{A_j})).} 
In the multiplicative version (i.e. for matrices $M_j$) 
it is stated by P.Deligne\footnote{It seems that the author of the present 
paper 
was the first to state the problem in the additive version} and C.Simpson 
was the first to obtain a significant result towards its resolution, see 
\cite{Si}. 

The problem is formulated without using the notions of fuchsian 
system and monodromy operator, yet they explain the interest in the 
problem. The multiplicative version is more important because   
the monodromy operators are invariant under the  
changes $X\mapsto W(t)X$ (where $W$ depends meromorphically on $t$ 
and det$W\not\equiv 0$) while the matrices-residua are 
not. In the multiplicative version the problem admits the interpretation: 
{\em for which $(p+1)$-tuples of local monodromies $c_j$ does there exist 
an irreducible monodromy group with such local monodromies ?}

The paper is structured as follows. In the next section we recall the basic 
results announced in \cite{Ko1} and proved in \cite{Ko2}. In 
Section~\ref{rigcase} we define the case of rigid $(p+1)$-tuples. In 
Section~\ref{examples} we give some new examples of existence of rigid 
$(p+1)$-tuples of diagonalizable matrices. In Section~\ref{examplesbis} we 
describe all rigid 
$(p+1)$-tuples of such matrices in which the multiplicities of all 
eigenvalues of one of the 
matrices are $\leq 2$. In Section~\ref{nondiagonal} we explain how the 
examples from the previous two sections give rise to other examples in which 
the matrices are not necessarily diagonalizable.

\section{The Deligne-Simpson problem for generic eigenvalues
\protect\label{DSPgen}}

\begin{defi}{\em A {\em Jordan normal form (JNF) of size $n$} is a family 
$J^n=\{ b_{i,l}\}$ ($i\in I_l$, $I_l=\{ 1,\ldots ,s_l\}$, $l\in L$) of 
positive integers $b_{i,l}$ 
whose sum is $n$. Here $L$ is the set of indices of the eigenvalues 
$\lambda _l$ (all distinct) and 
$I_l$ is the set of indices of Jordan blocks with eigenvalue $\lambda _l$;  
$b_{i,l}$ is the 
size of the $i$-th block with this eigenvalue. We assume that for each 
$l$ fixed one has $b_{1,l}\geq \ldots \geq b_{s_l,l}$. An $n\times n$-matrix 
$Y$ has the JNF $J^n$ (notation: $J(Y)=J^n$) if to its distinct 
eigenvalues $\lambda _l$, $l\in L$, there belong Jordan blocks of sizes 
$b_{i,l}$.}
\end{defi}

In what follows we presume the necessary condition $\prod \det (C_j)=1$ (resp. 
$\sum$Tr$(c_j)=0$) to hold. In terms of the eigenvalues $\sigma _{k,j}$  
(resp. $\lambda _{k,j}$) of the matrices from $C_j$ (resp. $c_j$) repeated 
with their multiplicities, this condition reads     

\[ \prod _{k=1}^n\prod _{j=1}^{p+1}\sigma _{k,j}=1~,~  
{\rm resp.~} \sum _{k=1}^n\sum _{j=1}^{p+1}\lambda _{k,j}=0~.\]  

\begin{defi} {\em An equality of the form  

\[ \prod _{j=1}^{p+1}\prod _{k\in \Phi _j}\sigma _{k,j}=1~,~{\rm resp.~} 
\sum _{j=1}^{p+1}\sum _{k\in \Phi _j}\lambda _{k,j}=0~,\]
is called a 
{\em non-genericity relation};  
the sets $\Phi _j$ contain one and the same number $\kappa$ ($1<\kappa <n$) 
of indices  
for all $j$. Eigenvalues satisfying none of these relations are called 
{\em generic}. Reducible  
$(p+1)$-tuples exist only for non-generic eigenvalues. Indeed, a reducible 
$(p+1)$-tuple can be conjugated to a block upper-triangular form and the 
eigenvalues of the restriction of the $(p+1)$-tuple to each diagonal block 
must satisfy condition (\ref{M_j}) or 
(\ref{A_j}) which is a  non-genericity relation.}
\end{defi}

For a conjugacy class $C$ in $GL(n,{\bf C})$ or $gl(n,{\bf C})$ denote by 
$d(C)$ its dimension. Remind that $d(C)$ is always even. 
For a matrix $Y$ from $C$ set 
$r(C):=\min _{\lambda \in {\bf C}}{\rm rank}(Y-\lambda I)$. The integer 
$n-r(C)$ equals the maximal number of Jordan blocks of $J(Y)$ with one and the 
same eigenvalue. 

Set $d_j:=d(C_j)$ (resp. $d(c_j)$), $r_j:=r(C_j)$ 
(resp. $r(c_j)$). The quantities 
$r(C)$ and $d(C)$ depend only on the JNF $J(Y)=J^n$, not 
on the eigenvalues and we write sometimes $r(J^n)$ and $d(J^n)$. 

\begin{prop}\label{d_jr_j}
(C. Simpson, see \cite{Si}.) The 
following couple of inequalities is a necessary condition for the existence 
of irreducible $(p+1)$-tuples of matrices $M_j$ satisfying (\ref{M_j}):

\[ d_1+\ldots +d_{p+1}\geq 2n^2-2~~~~~(\alpha _n)\] 
\[ {\rm for~all~}j,~r_1+\ldots +\hat{r}_j+\ldots +r_{p+1}\geq n~~~~~
(\beta _n)~~~.\]
\end{prop}

\begin{defi}{\em Denote by $\{ J_j^n\}$  
a $(p+1)$-tuple of JNFs, $j=1$,$\ldots$, $p+1$. 
We say that the DSP is {\em solvable} (resp. is {\em weakly solvable}) for 
a given $\{ J_j^n\}$ and given eigenvalues if there exists an 
irreducible $(p+1)$-tuple (resp. a $(p+1)$-tuple with a trivial centralizer) 
of matrices $M_j$ satisfying (\ref{M_j}) or of matrices $A_j$ satisfying 
(\ref{A_j}), with $J(M_j)=J_j^n$ or $J(A_j)=J_j^n$ and with the given 
eigenvalues. By definition, the DSP is solvable and weakly solvable for 
$n=1$.}
\end{defi}

\begin{tm}\label{2n}
The DSP is solvable for conjugacy classes $C_j$ or $c_j$ 
with generic eigenvalues and satisfying the condition 

\[ (r_1+\ldots +r_{p+1})\geq 2n~~~~~~~~~~~~~~~~(\omega _n)~~~.\] 
\end{tm}

For a given $\{ J_j^n\}$ with $n>1$, which satisfies conditions 
$(\alpha _n)$ and $(\beta _n)$ and doesn't satisfy condition 
$(\omega _n)$ set $n_1=r_1+\ldots +r_{p+1}-n$. Hence, $n_1<n$ and 
$n-n_1\leq n-r_j$. Define 
the $(p+1)$-tuple $\{ J_j^{n_1}\}$ as follows: to obtain the JNF 
$J_j^{n_1}$ 
from $J_j^n$ one finds one of the eigenvalues of $J_j^n$ with 
greatest number $n-r_j$ of Jordan blocks, then decreases  
by 1 the sizes of the $n-n_1$ {\em smallest} Jordan blocks with this 
eigenvalue and deletes the Jordan blocks of size 0. We use the notation 
$\Psi :\{ J_j^n\} \mapsto \{ J_j^{n_1}\}$. 

\begin{tm}\label{generic}
Let $n>1$. The DSP is solvable for the conjugacy classes $C_j$ or 
$c_j$ (with generic eigenvalues,  
defining the JNFs $J_j^n$ and satisfying conditions $(\alpha _n)$ and 
$(\beta _n)$) if and only if either $\{ J_j^n\}$ satisfies 
condition $(\omega _n)$ or the construction 
$\Psi :\{ J_j^n\}\mapsto \{ J_j^{n_1}\}$ iterated as long as it is defined 
stops at a $(p+1)$-tuple $\{ J_j^{n'}\}$ either with $n'=1$ or satisfying 
condition $(\omega _{n'})$.
\end{tm}

In the case of diagonalizable matrices $M_j$ or $A_j$ the JNF of $M_j$ or 
$A_j$ is completely defined by the {\em multiplicity vector (MV)} of its 
eigenvalues. This is a vector $\Lambda _j^n$ 
with positive integer components equal to the multiplicities of the 
eigenvalues of $M_j$ or $A_j$; hence, their sum 
is $n$. 

\begin{rem}\label{rd}
{\em Set $\Lambda _j^n=(m_{1,j},\ldots ,m_{i_j,j})$ where 
$m_{1,j}\geq \ldots \geq m_{i_j,j}$. Hence, one has  
$r_j=m_{2,j}+\ldots +m_{i_j,j}$ and $d_j=n^2-\sum _{i=1}^{i_j}(m_{i,j})^2$. 
In particular, the MV with greatest value of $d_j$ is $(1,\ldots ,1)$, with 
$d_j=n^2-n$.}
\end{rem}

Call {\em polymultiplicity vector (PMV)} the $(p+1)$-tuple of 
MVs $\Lambda ^n:=(\Lambda _1^n,\ldots ,\Lambda _{p+1}^n)$. 

\begin{rem}\label{diagcase}
{\em In the particular case of diagonalizable matrices $M_j$ or $A_j$ the 
mapping $\Psi$ is defined by the following rule (to be checked directly):}

The MV $\Lambda _j^{n_1}$ defining the JNF $J_j^{n_1}$ equals 
$(m_{1,j}-n+n_1,m_{2,j},m_{3,j},\ldots ,m_{i_j,j})$.
\end{rem}   

\section{The case of rigid $(p+1)$-tuples\protect\label{rigcase}}

\begin{defi}
{\em The case when $d_1+\ldots +d_{p+1}=2n^2-2$ is called {\em rigid}. 
Such $(p+1)$-tuples of matrices $A_j$ satisfying $(\ref{A_j})$ or of matrices 
$M_j$ satisfying $(\ref{M_j})$ or of JNFs or of PMVs are also called 
{\em rigid}.}
\end{defi}
 
A priori, if in the rigid case for a certain $(p+1)$-tuple of conjugacy 
classes 
the DSP is solvable, then up to conjugacy it has only finitely many 
solutions. 

\begin{prop}\label{SiKa}(see \cite{Si} and \cite{Ka})
If for a given $(p+1)$-tuple of conjugacy classes 
$C_j\subset GL(n,{\bf C})$ with generic eigenvalues and with 
$d_1+\ldots +d_{p+1}=2n^2-2$ the DSP is solvable for matrices $M_j$, then 
its solution is unique up to conjugacy.
\end{prop}

\begin{prop}\label{rigidA_j}
Suppose that for a given $(p+1)$-tuple of conjugacy classes 
$c_j\subset gl(n,{\bf C})$ with generic eigenvalues and with 
$d_1+\ldots +d_{p+1}=2n^2-2$ the DSP is solvable for matrices $A_j$. Then 
its solution is unique up to conjugacy.
\end{prop}

The proposition is proved at the end of the section.

\begin{rem}
{\em Rigid 
representations in the multiplicative case are studied in \cite{Ka} where an 
algorithm is given which 
tells whether the DSP is solvable for given conjugacy 
classes $C_j$ and the construction of rigid $(p+1)$-tuples of matrices 
$M_j\in C_j$ is explained. The algorithm of Katz is based on a middle 
convolution 
functor in the category of pervers sheaves. The same functor is defined in a 
purely algebraic way in \cite{DR}. The algorithm in \cite{Ka} also 
results in the construction 
$\Psi: \{ J_j^n\}\mapsto \{ J_j^{n_1}\}$ but in the case of rigid 
representations one never encounters $(p+1)$-tuples $\{ J_j^{n'}\}$ 
satisfying condition $(\omega _{n'})$. In fact, there holds the following 
lemma (see \cite{Ko1} and \cite{Ko2}):}
\end{rem}

\begin{lm}\label{invariant}
The quantity $2n^2-\sum _{j=1}^{p+1}d_j$ is invariant for the 
construction $\Psi :\{ J_j^n\}\mapsto \{ J_j^{n_1}\}$.
\end{lm}

The lemma implies that it is sufficient to check that condition $(\alpha _n)$ 
holds not for $\{ J_j^n\}$ (see Theorem~\ref{generic}) but for 
$\{ J_j^{n'}\}$. If $n'>1$ and condition $(\omega _{n'})$ holds, then for 
generic eigenvalues the DSP is solvable for the JNFs $J_j^{n'}$, see 
\cite{Ko1}, hence, the necessary condition $(\alpha _{n'})$ holds -- it is a 
strict inequality. If $n'=1$, 
then condition $(\alpha _{n'})$ is an equality (this 
is the rigid case). Hence, in both cases condition $(\alpha _{n'})$ holds and 
a posteriori one knows that in fact it is not necessary to check it. 
 
{\em Proof of Proposition~\ref{rigidA_j}:}

$1^0$. One can assume that for every $j$ there is no non-zero integer 
difference between two eigenvalues of the matrix $A_j$ (otherwise this can be 
achieved by a multiplication of the matrices by $c\in {\bf C}^*$). Hence, 

1) the monodromy 
operators $M_j$ of a fuchsian system with residua $A_j$ equal up to 
conjugacy $\exp (2\pi iA_j)$, see Remark~\ref{2pii};

2) the eigenvalues of the matrices $M_j$ are generic.

By Proposition~\ref{SiKa}, the $(p+1)$-tuple of matrices $M_j$ is 
unique up to conjugacy. Indeed, denote by $C_j$ the conjugacy class of $M_j$. 
Then $d(C_j)=d_j$ (see 1)) and $d(C_1)+\ldots +d(C_{p+1})=2n^2-2$. 

$2^0$. Suppose that there are at least two $(p+1)$-tuples of matrices 
$A_j\in c_j$ (denoted by $A_j^1$, $A_j^2$) 
non conjugate to one another which are solutions to the DSP. Denote by 
($F_1$), ($F_2$) two fuchsian systems with residua equal respectively to 
$A_j^1$, $A_j^2$ and with one and the same poles. Then these systems have one 
and the same monodromy group, see $1^0$. Hence, there exists a meromorphic 
change $X\mapsto W(t)X$ bringing ($F_1$) to ($F_2$). (The fact that $W$ can 
have a priori at most poles as singularities follows from the regularity 
of ($F_1$) and ($F_2$)).

$3^0$. For $t\neq a_j$, $j=1,\ldots ,p+1$, the matrix $W$ is holomorphic and 
holomorphically invertible. Indeed, it equals $X_2(X_1)^{-1}$ where $X_i$ is 
some fundamental solution to ($F_i$). Prove that $W$ has no pole at $a_j$. 

Suppose it has. Set $W=\sum _{k=-l}^{\infty }W_k(t-a_j)^k$, $l\in {\bf N}^*$. 
In a neighbourhood of $a_j$ one has 

\[ (F_i):{\rm d}X/{\rm d}t=(A_j^i/(t-a_j)+O(1))X~,~i=1,2.\]
Then one has 

\[ -W^{-1}{\rm d}W/{\rm d}t+W^{-1}(A_j^1/(t-a_j)+O(1))W=
(A_j^2/(t-a_j)+O(1))~,~{\rm i.e.}\] 

\[ -{\rm d}W/{\rm d}t+(A_j^1/(t-a_j)+O(1))W=
W(A_j^2/(t-a_j)+O(1))~~{\rm and, hence,}\]

\[ -lW_{-l}+A_j^1W_{-l}-W_{-l}A_j^2=0~.\] 
This implies that $W_{-l}=0$, i.e. $W$ has no pole at $a_j$. (Indeed, the 
eigenvalues of the linear operator $(.)\mapsto -l(.)+A_j^1(.)-(.)A_j^2$ 
acting on $gl(n,{\bf C})$ are of the form $\eta =-l+\lambda _1-\lambda _2$ 
where 
$\lambda _i$ is eigenvalue of $A_j^i$. Their set is one of the same for 
$i=1,2$ and by $1^0$, one has $\eta \neq 0$.)

But then $W$ is holomorphic on ${\bf C}P^1$, hence, 
constant, i.e. $W\in GL(n,{\bf C})$ which means that the two 
$(p+1)$-tuples (of matrices $A_j^1$ and $A_j^2$) are conjugate.

The proposition is proved.

\section{Some series of rigid representations\protect\label{examples}}

In this section we list several series of rigid representations with 
diagonalizable matrices $M_j$ or $A_j$ by means 
of their PMVs. Their existence follows from Theorem~\ref{generic} and 
Remark~\ref{diagcase}. (The eigenvalues are presumed generic.) In  
Section~\ref{nondiagonal} we explain how to deduce from their existence the 
one of other rigid series with generic eigenvalues in which at least one of 
the matrices $M_j$ or $A_j$ is not diagonalizable.

For $p=2$ we define several series of PMVs. We avoid the letters $A$ and $M$ 
which denote already matrices and the notation should not be mixed up with 
similar notation for singularities or Lie algebras:\\ 

$\begin{array}{llll}W_k:&(k,k,k+1),&(k,k,k+1),&(k,k,k+1)\\
B_k:&(k,k,k-1),&(k,k,k-1),&(k,k,k-1)\\
C_k:&(k,k,k),&(k,k,k),&(k,k+1,k-1)\\
D_k:&(k,k,k,k+1),&(k,k,k,k+1),&(2k,2k+1)\\
E_k:&(k,k,k,k-1),&(k,k,k,k-1),&(2k,2k-1)\\
F_k:&(k,k,k,k),&(k,k,k,k),&(2k+1,2k-1)\\
\Phi _k:&(k,k,k+1,k-1),&(k,k,k,k),&(2k,2k)\\
G_k:&(k,k,k+1,k+1),&(k,k,k+1,k+1),&(2k+1,2k+1)\\
H_k:&(k,k,k,k,k,k+1),&(3k,3k+1),&(2k,2k,2k+1)\\
I_k:&(k,k,k,k,k,k-1),&(3k,3k-1),&(2k,2k,2k-1)\\
J_k:&(k,k,k,k,k,k),&(3k+1,3k-1),&(2k,2k,2k)\\
K_k:&(k,k,k,k,k,k),&(3k,3k),&(2k,2k+1,2k-1)\\
L_k:&(k,k,k,k,k+1,k-1),&(3k,3k),&(2k,2k,2k)\\
V_k:&(k,k,k,k,k+1,k+1),&(3k+1,3k+1),&(2k,2k+1,2k+1)\\
N_k:&(k,k,k,k+1,k+1,k+1),&(3k+1,3k+2),&(2k+1,2k+1,2k+1)\\
P_k:&(k,k,k,k,k-1,k-1),&(3k-1,3k-1),&(2k,2k-1,2k-1)\end{array}$\\ 

Here $k\in {\bf N}$ or $k\in {\bf N}^+$ according to the case. Each of these 
PMVs satisfies Conditions $(\alpha _n)$ and $(\beta _n)$ (to be checked 
directly). Moreover, $(\alpha _n)$ is equality everywhere. 

The series $W_k$, $B_k$ and $C_k$ were discovered by O.Gleizer (see 
\cite{Gl}). We don't use his result but deduce their existence from 
Theorem~\ref{generic} and Remark~\ref{diagcase} (partly because we need to 
prove the existence of 
rigid triples from other 
series as well and partly because he claims in \cite{Gl} the non-existence 
of the rigid series 

\[ OG_k:~(2,\ldots ,2,1,1,1)~,~(2,\ldots ,2,1)~,~(2k-1,1,1)\]  
which contradicts Theorem~\ref{generic}; we deduce the 
existence of this series at the end of the section.). 

To prove the existence of these rigid series it suffices to explicit the 
sequence of PMVs $\Lambda ^n$, $\Lambda ^{n_1}$, $\ldots$, $\Lambda ^{n_s}$ 
occurring when the construction $\Psi$ from Section~\ref{DSPgen} is iterated, 
see Theorem~\ref{generic} and Remark~\ref{diagcase}; we set $n_s=n'$. For 
$\Lambda ^n=W_k$ this sequence equals 
$W_k$, $B_k$, $W_{k-1}$, $B_{k-1}$, $\ldots$, $W_0$. Write it symbolically 
in the form 

\[ W_k\rightarrow B_k\rightarrow W_{k-1}\rightarrow B_{k-1}\rightarrow 
\ldots \rightarrow W_0\]
All requirements of Theorem~\ref{generic} and Remark~\ref{diagcase} are met 
which implies the 
existence of irreducible triples with PMV $W_k$ (and $B_k$ as well if one 
deletes the first term of the sequence). 

One finds by analogy (for $\Lambda ^n=C_k$) the sequence

\[ C_k\rightarrow B_k\rightarrow W_{k-1}\rightarrow B_{k-1}\rightarrow 
\ldots \rightarrow W_0\]
which differs from the previous one only in its first term. Hence, there 
exist irreducible triples with $\Lambda ^n=C_k$. In the same way one obtains 
the sequences
 
\[ D_k~{\rm or}~F_k~{\rm or}~\Phi _k\rightarrow E_k\rightarrow G_{k-1}
\rightarrow D_{k-1}
\rightarrow E_{k-1}\rightarrow G_{k-2}\rightarrow \ldots \rightarrow G_0
\rightarrow D_0\]

\[ H_k~{\rm or}~J_k~{\rm or}~K_k~{\rm or}~L_k\rightarrow I_k\rightarrow P_k
\rightarrow N_{k-1}\rightarrow 
V_{k-1}\rightarrow H_{k-1}\rightarrow I_{k-1}
\rightarrow \ldots \rightarrow 
H_0\]
from which one deduces the existence of irreducible triples with 
$\Lambda ^n$ equal to any of the other PMVs listed above.

For $p=3$ we define two series:

$R_k$: $(k,k)$, $(k,k)$, $(k,k)$, $(k+1,k-1)$;

$S_k$: $(k+1,k)$, $(k+1,k)$, $(k+1,k)$, 
$(k+1,k)$.

The corresponding sequence equals 

\[ S_k~{\rm or}~R_k\rightarrow S_{k-1}\rightarrow S_{k-2}\rightarrow \ldots 
\rightarrow S_0\]

For $p=4$ we define the series

\[ T_k~:~(2k+1,2k-1)~,~(3k,k)~,~(3k,k)~,~(3k,k)~,~(3k,k)~.\] 

The PMV $\Lambda ^{n_1}$ equals $(2k-1)$, $(k,k-1)$, $(k,k-1)$, $(k,k-1)$, 
$(k,k-1)$. This means that the matrix $A_1$ must be scalar and the PMV of 
the other four matrices equals $S_{k-1}$. Thus, the existence of irreducible 
quintuples follows from the existence of irreducible quadruples with PMV 
$S_{k-1}$. 

Finally, we recall the existence of other four series discovered by 
C.Simpson (the first three, see \cite{Si}) and by O.Gleizer (see \cite{Gl}):

\[ \begin{array}{lllll}HG_n~:&(n-1,1)&(1,\ldots ,1)&(1,\ldots ,1)&
{\rm hypergeometric}\\
OF_n~:&((n+1)/2, (n-1)/2)&((n-1)/2,(n-1)/2,1)&(1,\ldots ,1)&{\rm odd~family}\\
EF_n~:&(n/2,n/2)&(n/2,(n-2)/2,1)&(1,\ldots ,1)&{\rm even~family}\\
FF_n~:&(2,1,\ldots ,1)&(2,\ldots,2,1,\ldots ,1)& 
(n-2,2)&{\rm finite~family~,~}\\&&(n-4~{\rm times~}2)&&n=5,6,7,8\end{array}\]

For the series $OG_k$ defined above one obtains the sequence

\[ OG_k\rightarrow OG_{k-1}\rightarrow \ldots \rightarrow OG_1\rightarrow 
HG_2\rightarrow HG_1\]
Note that $OG_1=HG_3$.

The existence of the series

\[ [n-1,1]~:(n-1,1)~,\ldots ,~(n-1,1)~~(n+1~{\rm times~})\]
follows from $[n-1,1]\rightarrow (1),\ldots ,(1)$.

\begin{rem}\label{solutions}
{\em In the series $W_k - P_k$ the multiplicities of the eigenvalues 
are equal to (or differ by no more than 2 from) $n/s_1$, $n/s_2$, $n/s_3$ 
where $(s_1,s_2,s_3)\in ({\bf N}^*)^3$ is a solution to the equation}

\[ 1/s_1+1/s_2+1/s_3=1\]
{\em (these solutions are (3,3,3), (4,4,2) and (6,3,2) up to permutation). 
One can 
consider the series $OF_n$ and $EF_n$ (resp. $HG_n$) 
as corresponding to the ``generalized'' 
solution $(2,2,\infty )$ (resp. $(1,\infty ,\infty )$) of the above equation.}
\end{rem}

\begin{rem}\label{3series} 
{\em C.Simpson has shown 
in \cite{Si} that the three series $OF_n$, $EF_n$ and $HG_n$ include all 
rigid triples of 
diagonalizable matrices $M_j$ in which one of them has distinct eigenvalues. 
Hence, this is the case of matrices $A_j$ as well because the criterium for 
existence of irreducible $(p+1)$-tuples (i.e. Theorem~\ref{generic}) is the 
same in the additive and in the multiplicative situation.}
\end{rem} 

\section{Rigid representations with an upper bound on the multiplicities 
of the eigenvalues of the first matrix\protect\label{examplesbis}}

\subsection{Formulation of the problem}

In the present section we consider the problem: 

{\em Give the complete list of PMVs for which there exist rigid 
irreducible $(p+1)$-tuples of diagonalizable matrices $M_j$ satisfying 
(\ref{M_j}) (resp. of diagonalizable matrices $A_j$ satisfying 
(\ref{A_j})), with generic eigenvalues, in which the multiplicities of all 
eigenvalues of $M_1$ (resp. of $A_1$) are $\leq u$ for some $u\in {\bf N}^*$.}

We solve the problem for $u=2$. In what follows we set $m_{1,1}=u=2$. 
(If $m_{1,1}=1$, then $u=1$ and in this case the answer to the 
problem is given by Remark~\ref{3series}.) The techniques can be used to 
solve the problem for any given $u$. We assume that no MV equals ($n$) 
in which case the corresponding matrix $A_j$ or $M_j$ must be scalar. We also 
assume that no MV is of the form $(1,\ldots ,1)$ (see Remark~\ref{3series}). 

\begin{rem}
{\em The cases $u=1$ 
and $u=2$ are exceptional in the following sense -- whenever one finds 
a rigid PMV satisfying condition $(\beta _j)$, there exist rigid 
$(p+1)$-tuples of diagonalizable matrices with this PMV. (For $u=3$ this is 
not true, see Example~\ref{exu=3}.) More generally, there holds}
\end{rem}

\begin{tm}\label{uleq2}
If $u\leq 2$, then conditions $(\alpha _n)$ and $(\beta _n)$ are necessary and 
sufficient for the existence for generic eigenvalues 
of irreducible $(p+1)$-tuples of matrices $M_j$ satisfying (\ref{M_j}) or 
of matrices $A_j$ satisfying (\ref{A_j}).
\end{tm}

The theorem is proved in Section~\ref{pruleq2}. It 
generalizes Simpson's result from \cite{Si}: {\em if 
one of the matrices $M_j$ has distinct eigenvalues, then for generic 
eigenvalues conditions $(\alpha _n)$ and $(\beta _n)$ are necessary and 
sufficient for the existence  
of irreducible $(p+1)$-tuples of matrices $M_j$ satisfying (\ref{M_j}).} 
In the above theorem condition $(\alpha _n)$ is not presumed to be an 
equality, i.e. the theorem does not consider only the rigid case. 

\begin{ex}\label{exu=3}
{\em For $p=2$, $u=3$, $n=6m+3$, $m\in {\bf N}^*$ the PMV 
$(3,\ldots ,3,2,1,\ldots ,1)$ ( $m$ times 3, $3m+1$ units), 
$(3m+1,3m+1,1)$, $(3m+1,3m+1,1)$ is 
rigid and satisfies condition $(\beta _n)$ but the PMV obtained from it after 
applying the construction $\Psi$ (one has $n_1=n-2$) 
does not satisfy condition $(\beta _{n-2})$.}
\end{ex}

\subsection{The results\protect\label{results}}

The basic result is contained in Theorems \ref{p=3}, \ref{pleq3} and 
\ref{p=2}. In the next subsection we explain the method of proof.  

\begin{tm}\label{p=3}
If $u=2$, $p=3$ and $\Lambda _4^n=(n-1,1)$, then 

1) one has $d_1+\ldots +d_4\geq 2n^2-2$ 
in all cases except in 

{\bf Case} $\Omega$: $n$ is even, 
$r_2=n/2$, $r_3=n/2-1$, 
$\Lambda _1^n=(2,\ldots ,2)$, $\Lambda _2^n=(n/2,n/2)$, 
$\Lambda _3^n=(n/2+1,n/2-1)$; 
in {\bf Case} $\Omega$ one 
has $d_1+d_2+d_3+d_4=2n^2-4$;    

2) the only PMVs of rigid 
quadruples for $n$ even are 

\[ \begin{array}{lllll}\Xi _n~:&(2,\ldots ,2)&(n/2,n/2)&(n/2,n/2)&(n-1,1)\\ 
\Theta _n~:&(2,\ldots ,2,1,1)&(n/2,n/2)&(n/2+1,n/2-1)&(n-1,1)\\
\Psi _6~:&(2,2,2)&(3,3)&(4,1,1)&(5,1)\end{array}\] 
and the only ones for $n$ odd are 

\[ \begin{array}{lllll}\Pi _n~:&(2,\ldots ,2,1)&((n+1)/2,(n-1)/2)&
((n+1)/2,(n-1)/2)&(n-1,1)\\
\Delta _n~:&(2,\ldots ,2,1)&(2,\ldots ,2,1)&(n-1,1)&(n-1,1)\end{array}\]
\end{tm}

The theorem is proved in Subsection \ref{p=3pr}. 

\begin{tm}\label{pleq3}
If $u=2$, then for a rigid 
$(p+1)$-tuple one has $p\leq 3$. If $p=3$, then one of the MVs of a rigid 
quadruple equals $(n-1,1)$. 
\end{tm}

The theorem is proved in Subsection \ref{pleq3pr}.

\begin{tm}\label{p=2}
If $u=2$ and $p=2$, then with 
the exception of finitely many cases with $n\leq 21$ the 
only PMVs for which there exist rigid triples are the following ones:

For $n$ even: 
\[ \begin{array}{llllll}1a)&\Gamma ^1_n~:&(2,\ldots ,2)&
(2,\ldots ,2,1,1,1,1,1,1)&(n-2,2)\\ 
1b)&\Gamma ^2_n~:&(2,\ldots ,2,1,1)&(2,\ldots ,2,1,1,1,1)&(n-2,2)\\ 
1c)&\Gamma ^3_n~:&(2,\ldots ,2,1,1)&(2,\ldots ,2,1,1)&(n-2,1,1)\\ 
1d)&\Gamma ^4_n~:&(2,\ldots ,2)&(2,\ldots ,2,1,1,1,1)&(n-2,1,1)\\
1e)&Y^1_n~:&(2,\ldots ,2,1,1,1,1)&(m,m,2)&(m+1,m+1)\\
1f)&Y^2_n~:&(2,\ldots ,2,1,1)&(m,m,1,1)&(m+1,m+1)\\
1g)&Y^3_n~:&(2,\ldots ,2,1,1,1,1)&(m+2,m,1,1)&(m+2,m+2)\\
1h)&Y^4_n~:&(2,\ldots ,2,1,1,1,1,1,1)&(m+2,m,2)&(m+2,m+2)\\
1i)&Y^5_n~:&(2,\ldots ,2,1,1)&(m+1,m,1)&(m+1,m,1)\\
1j)&Y^6_n~:&(2,\ldots ,2,1,1,1,1)&(m,m,1,1)&(m+2,m)\\
1k)&Y^7_n~:&(2,\ldots ,2,1,1,1,1,1,1)&(m,m,2)&(m+2,m)
\end{array}\] 

For $n$ odd: 
\[ \begin{array}{lllll}2a)&X^1_n~:&(2,\ldots ,2,1,1,1,1,1)&
(2,\ldots ,2,1)&(n-2,2)\\
2b)&X^2_n~:&(2,\ldots ,2,1,1,1)&(2,\ldots ,2,1,1,1)&(n-2,2)\\
2c)&OG_n~:&(2,\ldots ,2,1)&(2,\ldots ,2,1,1,1)&(n-2,1,1)\\ 
2d)&Z^1_n~:&(2,\ldots ,2,1)&(m,m,1)&(m,m,1)\\
2e)&Z^2_n~:&(2,\ldots ,2,1,1,1,1,1)&(m,m-1,2)&(m+1,m)\\
2f)&Z^3_n~:&(2,\ldots ,2,1,1,1)&(m,m-1,1,1)&(m+1,m)\\
2g)&Z^4_n~:&(2,\ldots ,2,1,1,1)&(m,m,1)&(m+1,m-1,1)
\end{array}\]
where $m\in {\bf N}$ or $m\in {\bf N}^*$.
\end{tm}

The theorem is proved in Subsection \ref{p=2pr}. We do not explicit the 
exceptional cases with $n\leq 21$. The reader can do this by iterating 
the construction $\Psi$ from Section~\ref{DSPgen} backward. 

\subsection{The method of proof\protect\label{method}}

The method of proof consists in trying to minimize the 
quantities $d_j$ for $r_j$ fixed. Denote these minimal possible values of 
$d_j$ by $d_j'$ 
and the PMVs realizing these minimal values by ${\Lambda '}^n$. (A posteriori 
they turn out to be unique up to permutation of the components of their MVs.) 

The PMVs ${\Lambda '}^n$ in part of the cases turn out to 
be rigid and then we prove the existence of the corresponding $(p+1)$-tuples 
of matrices by means of Theorem~\ref{generic}. In another part of the cases 
one finds out that $d_1'+\ldots +d_{p+1}'>2n^2-2$, i.e. no rigid 
$(p+1)$-tuples exist for such quantities $r_j$. Finally, in the remaining part 
of the cases one has $d_1'+\ldots +d_{p+1}'<2n^2-2$ (i.e. no irreducible 
$(p+1)$-tuples of matrices exist for the PMVs ${\Lambda '}^n$) and one finds 
out how to change the PMVs in order to have 
$d_1'+\ldots +d_{p+1}'=2n^2-2$, without changing the quantities $r_j$; 
after this one proves the existence of 
rigid $(p+1)$-tuples from the new PMVs.

The following lemmas explain how 
this is done in more details. Recall that we denote by $u$ the component 
$m_{1,1}$ of $\Lambda _1^n$ and 
that $m_{1,1}\geq \ldots \geq m_{i_1,1}$. 

\begin{lm}\label{min1}
If $r_j\leq n/2$ is fixed, then $d_j$ is minimal if and only if 
$\Lambda _j^n=(n-r_j,r_j)$.
\end{lm}

{\em Proof:}

One has $d_j=n^2-(m_{1,j})^2-\sum _{k=2}^{i_j}(m_{k,j})^2$ where 
$m_{1,j}=(n-r_j)\geq n/2$, see 
Remark~\ref{rd}. The sum $\sum _{k=2}^{i_j}(m_{k,j})^2$ is maximal if and 
only if $i_j=2$, $m_{2,j}=r_j$.

\begin{defi}
{\em Recall that $\Lambda _j^n=(m_{1,j},\ldots ,m_{i_j,j})$, 
$m_{1,j}\geq \ldots \geq m_{i_j,j}$. If 
$m_{1,j}=\ldots =m_{\mu ,j}>m_{\mu +1,j}$, $\mu +1<i_j$, then the change 
$m_{\mu +1,j}\mapsto m_{\mu +1,j}+1$, $m_{i_j,j}\mapsto m_{i_j,j}-1$ is 
called a {\em passage}. Its inverse is called an {\em antipassage}. A 
passage preserves $r_j$ and decreases 
$d_j$ (to be checked directly). If after the change one has $m_{i_j,j}=0$, 
then one deletes the last component of $\Lambda _j^n$ and sets 
$i_j\mapsto i_j-1$.}
\end{defi}

\begin{lm}\label{min2}
If $r_j>n/2$ is fixed, then $d_j$ is minimal if and only if 
$\Lambda _j^n=(m,m,\ldots ,m,q)$ where $1\leq q\leq m=n-r_j$.
\end{lm}

{\em Proof:}

Perform passages as long as they are defined. 
No matter what the components $m_{i,j}$ are at the 
beginning, at the end one has $\Lambda _j^n=(m_{1,j},\ldots ,m_{1,j},q)$.

\begin{cor}\label{min3}
If $u=2$, then $d_1$ is minimal if 
and only if $\Lambda _1^n=(2,\ldots ,2)$ for $n$ even (and, hence, 
$d_1=n^2-2n$) or 
$\Lambda _1^n=(2,\ldots ,2,1)$ for $n$ odd and $d_1=n^2-2n+1$. 
\end{cor}

The corollary is direct.

\begin{rem}\label{twoMVs}
{\em Suppose that two of the MVs equal $(\alpha ,\beta )$, $(v,w)$ with 
$\alpha >\beta$, $v\geq w$, $\beta \geq w$ (hence, $\alpha \leq v$) and 
$\beta +1\leq n/2$. 
Hence, their quantities $d_j$ equal respectively $2\alpha \beta$, $2vw$. 
Their quantities $r_j$ equal respectively $\beta$, $w$. 

Change the two MVs to $(\alpha -1,\beta +1)$, $(v+1,w-1)$. Hence, their 
new quantities $r_j$ are $\beta +1$, $w-1$, i.e. their sum does not change. 
The new quantities $d_j$ are $\alpha \beta +\alpha -\beta -1$, $vw+w-v-1$, 
their sum changes by $\alpha -\beta +w-v-2=(\alpha -v)+(w-\beta )-2<0$, i.e. 
their sum decreases.}
\end{rem}

\begin{lm}\label{sumr}
If $u=2$, then one has $r_2+\ldots +r_{p+1}=n$ or $n+1$.
\end{lm}

Indeed, if $u=2$, then $r_1=n-2$. For rigid $(p+1)$-tuples condition 
$(\beta _n)$ holds while condition 
$(\omega _n)$ does not. This leaves only the two possible values ($n$ and 
$n+1$) for $r_2+\ldots +r_{p+1}$.

\subsection{Proof of Theorem~\protect\ref{p=3}\protect\label{p=3pr}} 

$1^0$. To prove the theorem we consider all cases in which for 
given quantities $r_j$ the corresponding quantities $d_j$ are minimal. 
They are given by Lemmas~\ref{min1}, \ref{min2} and 
Corollary~\ref{min3}. We prove that among these cases {\bf Case} $\Omega$ is 
the only one in which condition $(\alpha _n)$ does not hold. This is part 1) 
of the theorem. We also find all rigid cases among them (this is part 2)). 

$2^0$. There are two possible cases: $r_2+r_3=n-1$ or $n$ (Lemma~\ref{sumr}).  
 
{\bf Case 1)} $r_2+r_3=n$. 

{\bf Subcase 1.1)} $r_2=r_3=n/2$ (i.e. $n$ is even). 

One has $\Lambda _2^n=\Lambda _3^n=(n/2,n/2)$, $d_2=d_3=n^2/2$  
(Lemma~\ref{min1}), $d_4=2n-2$  
and $d_1\geq n^2-2n$ (Corollary~\ref{min3}). Hence, to have rigid 
quadruples the last inequality 
must be equality and we have the series $\Xi _n$. 

{\bf Subcase 1.2)} $r_2>r_3$. 

One has $r_2>n/2$, $r_3<n/2$ and by Lemmas~\ref{min2} and \ref{min1} 
$d_2$, $d_3$ are minimal if and only if 
$\Lambda _2^n=(m,\ldots ,m,s)$, $\Lambda _3^n=(n-m,m)$ where $n=lm+s$, 
$l\in {\bf N}$, $1\leq s\leq m$, $r_2=(l-1)m+s$, $r_3=m$.

Hence, $l\geq 2$ (otherwise $r_2\leq n/2$). One has 

\[ d_1\geq n^2-2n~,~d_2=l(l-1)m^2+2lms~,~d_3=2m((l-1)m+s)~{\rm and~}
d_4=2n-2~.\] 
Set $\Delta = d_1+d_2+d_3+d_4-(2n^2-2)$. Hence, 

\[ \Delta \geq 
-n^2-2n+2+l(l-1)m^2+2lms+2m((l-1)m+s)+2n-2=\]
\[ =-(ml+s)^2+l(l-1)m^2+2lms+2m((l-1)m+s)=\]
\[ =(l-2)m^2+2ms-s^2=(l-2)m^2+ms+s(m-s)>0~.\]
This means that rigid $(p+1)$-tuples with $r_2>n/2$, $r_3<n/2$ and 
$r_2+r_3=n$ do not exist. 

{\bf Case 2)} $r_2+r_3=n-1$. 

{\bf Subcase 2.1)} $r_2=r_3=(n-1)/2$ (i.e. $n$ is odd).
 
One has $d_2=d_3=(n^2-1)/2$, $d_4=2n-2$ and $d_1\geq n^2-2n+1$ with 
equality if and only if $\Lambda _1^n=(2,\ldots ,2,1)$ 
(Corollary~\ref{min3}). Hence, to have a rigid quadruple the last 
inequality must be equality and we have the series $\Pi _n$. 

{\bf Subcase 2.2)} $n$ is odd and $r_2>r_3$. 

One has $r_2>n/2$, $r_3<n/2$ and by Lemmas~\ref{min2} and \ref{min1} 
$d_2$, $d_3$ are minimal if and only if 
$\Lambda _2^n=(m,\ldots ,m,s)$, $\Lambda _3^n=(n-m+1,m-1)$ where $n=lm+s$, 
$l\in {\bf N}$, $1\leq s\leq m$, $r_2=(l-1)m+s$, $r_3=m-1$. 

Hence, $l\geq 2$, 
(otherwise $r_2<n/2$; note that $l=1$, $s=m$ is impossible because $n$ is odd) 
and $m>1$ (otherwise $A_3$ or $M_3$ must be scalar).
One has 

\[ \Delta \geq 
-n^2-2n+2+l(l-1)m^2+2lms+2(m-1)((l-1)m+s+1)+2n-2=\]
\[ =-(ml+s)^2+l(l-1)m^2+2lms+2(m-1)((l-1)m+s+1)=\]
\[ =(l-2)m^2+2ms-s^2-2(l-1)m-2s-2+2m=(l-2)m^2+ms+s(m-s)-2(l-2)m-2s-2=\]
\[ =(l-2)m(m-2)+s(m-s)+(m-2)s-2>0\]
for $m>2$ because either $l>2$ or $l=2$ and $m\geq s\geq 1$. Hence, there 
are no rigid quadruples in this case. If $m=2$, then 
$\Lambda _3=\Lambda _4=(n-1,1)$ -- this gives the series $\Delta _n$.

{\bf Subcase 2.3)} $n$ is even, $r_2>r_3$,   
$r_2>n/2$ and $r_3<n/2$.

Like in {\bf Subcase 2.2)} we show that no rigid quadruples 
exist (it is impossible to have $l=1$, $s=m$ because in this case $r_2=n/2$). 

{\bf Subcase 2.4)} $n$ is even, $r_2=n/2$, $r_3=n/2-1$ and 
$\Lambda _1^n=(2,\ldots ,2)$, $\Lambda _2^n=(n/2,n/2)$, 
$\Lambda _3^n=(n/2+1,n/2-1)$. 

One has $d_1+d_2+d_3+d_4=2n^2-4$. 
This is precisely {\bf Case} $\Omega$. In 
this case to have an irreducible representation 
one cannot choose for all three matrices $A_1$, $A_2$, $A_3$ (or 
$M_1$, $M_2$, $M_3$) the Jordan normal forms defined by the 
MVs $\Lambda _1^n$, $\Lambda _2^n$, $\Lambda _3^n$.

All conjugacy classes are even-dimensional. To have a rigid quadruple  
one has to choose only for one of the indices $j=1,2,3$ a conjugacy class of 
dimension $\rho$ next after the minimal one $\rho _{{\rm min}}=d_j$ and 
one must have $\rho =\rho _{{\rm min}}+2$ (because $d_1+\ldots +d_4$ has to 
increase by 2). 

For $n\geq 8$ this can be done 
only for $j=1$ and this gives the series $\Theta _n$. For $n=6$ one can choose 
$j=3$ as well (but not $j=2$) and this gives the case $\Psi _6$. For $n=4$ 
the only possibility is $(2,1,1)$, $(2,2)$, $(3,1)$, $(3,1)$ which is the case 
$\Theta _4$.   
 
$3^0$. Prove that rigid quadruples from the five cases $\Xi _n$, $\Theta _n$, 
$\Psi _6$, $\Pi _n$ and $\Delta _n$ really exist. 
Use the notation from Section~\ref{examples}. One has 

\[ \Xi _{n+1}\rightarrow \Pi _n~{\rm and~}
\Pi_n\rightarrow \Pi _{n-2}\rightarrow \ldots \rightarrow 
\Pi _3=[2,1]~,\]
this proves the existence of the rigid series $\Xi _n$ and $\Pi _n$. One also 
has 

\[ \Theta _n\rightarrow \Theta _{n-2}\rightarrow \ldots \rightarrow \Theta _4
\rightarrow HG_2~{\rm and~}\Psi _6\rightarrow \Theta _4\rightarrow HG_2~.\] 
This proves the existence of the series $\Theta _n$ and of $\Psi _6$. The one 
of the series $\Delta _n$ follows from $\Delta _n\rightarrow \Delta _{n-2}$ 
and $\Delta _3=[2,1]$.

The theorem is proved.

\subsection{Proof of Theorem~\protect\ref{pleq3}\protect\label{pleq3pr}}

$1^0$. Recall that the change of two MVs $(n-r_j,r_j)$, $(n-r_i,r_i)$  
to $(n-r_j-1,r_j+1)$, $(n-r_i+1,r_i-1)$ (provided that $r_i\leq r_j\leq n/2$ 
and $r_j+1\leq n/2$) does 
not change the sum $r_j+r_i$ and decreases the sum $d_j+d_i$, see 
Remark~\ref{twoMVs}. In what follows when such a change is performed and 
after it a MV becomes equal to $(n)$ we delete it because the corresponding 
matrix $A_j$ or $M_j$ must be scalar.

Remind that MVs like the above ones give the minimal value of $d_j$ when 
$r_j$ is fixed and $r_j\leq n/2$, see Lemma~\ref{min1}.

$2^0$. Consider only these $(p+1)$-tuples ($p\geq 4$) in which the MVs 
provide minimal 
possible values for $d_j$ when $r_j$ is fixed (see Lemmas~\ref{min1} and 
\ref{min2} and Corollary~\ref{min3}). For all of them we show that 
condition $(\alpha _n)$ holds and is a strong inequality. Hence, it is strong 
for all other possible MVs with these values of $r_j$, i.e. no 
rigid $(p+1)$-tuples exist for $p\geq 4$. 

As a result of suitably chosen changes of MVs like in $1^0$ one comes to the 
case $p=3$, $\Lambda _4^n=(n-1,1)$. In this case one has 
$d_1+\ldots +d_4<2n^2-2$ only in  
{\bf Case} $\Omega$, see Theorem~\ref{p=3}, when one has 
$d_1+\ldots +d_4=2n^2-4$. 

Hence, if starting with a $(p+1)$-tuple one comes as a result of such 
changes of MVs to the case 
$p=3$, $\Lambda _4^n=(n-1,1)$, but not to {\bf Case} $\Omega$, then the 
$(p+1)$-tuple is not rigid, see Remark~\ref{twoMVs}. 

$3^0$. So consider only the $(p+1)$-tuples which after a change like in $1^0$ 
become the quadruple from {\bf Case} $\Omega$. This means that either 
$p=4$ or $p=3$ (as a result of a change of MVs no more than one MV of the form 
$(n)$ can appear). 

We show in $4^0$ why the case $p=3$ needs not to be considered. 
If $p=4$, then there are only two possibilities:

\[ \begin{array}{llll}1)&\Lambda _1^n=(2,\ldots ,2)&
\Lambda _2^n=\Lambda _3^n=(n/2+1,n/2-1)&\Lambda _4^n=\Lambda _5^n=(n-1,1)\\
2)&\Lambda _1^n=(2,\ldots ,2)&
\Lambda _2^n=(n/2,n/2)~,~\Lambda _3^n=(n/2+2,n/2-2)&
\Lambda _4^n=\Lambda _5^n=(n-1,1)~.\end{array}\]  

One has respectively $d_1+\ldots +d_5=2n^2+2n-8$ and $2n^2+2n-12$. Hence, 
the first possibility never gives a rigid quintuple (one has $n\geq 4$). 
The second can give a rigid quintuple only for $n=5$, but $n$ must be even. 
Note that for $n=4$ the MV $\Lambda _3^n$ from 2) equals (4), so this is 
in fact a quadruple, not a quintuple.

$4^0$. If as a result of changes of MVs a $(p+1)$-tuple with $p\geq 4$ becomes 
first a quadruple different from the one of {\bf Case} $\Omega$ and then 
the one from {\bf Case} $\Omega$, then it cannot be rigid -- each change 
decreases $d_1+\ldots +d_{p+1}$ by at least 2 and in {\bf Case} $\Omega$ this 
sum equals $2n^2-4$.  

The theorem is proved.

\subsection{Proof of Theorem \protect\ref{p=2}\protect\label{p=2pr}}

Set $\Delta = d_1+d_2+d_3-(2n^2-2)$. Irreducible (resp. rigid) triples 
can exist only for $\Delta \geq 0$ (resp. $\Delta =0$), see condition 
$(\alpha _n)$.  Like in the proof of Theorem~\ref{p=3} we consider all cases 
in which for given quantities $r_j$ the corresponding quantities $d_j$ are 
minimal, see Lemmas~\ref{min1}, \ref{min2} and 
Corollary~\ref{min3}. We assume that no MV equals $(1,\ldots ,1)$, see 
Remark~\ref{3series}.

{\bf Case 1)} $r_2+r_3=n$.

{\bf Subcase 1.1)} $r_2>n/2$, $r_3<n/2$. 

$1^0$. The quantities $d_2$ and $d_3$ are minimal 
if and only if one has $\Lambda _2^n=(m,\ldots ,m,s)$, $\Lambda _3^n=(n-m,m)$, 
$n=lm+s$, $l\in {\bf N}$, $1\leq s\leq m$, $r_2=(l-1)m+s$, $r_3=m$, 
see Lemmas~\ref{min1} and \ref{min2}. For such $\Lambda _2^n$, $\Lambda _3^n$ 
one has $d_2=l(l-1)m^2+2lms$, $d_3=2m((l-1)m+s)$. Hence, 

\[ \Delta \geq 
-n^2-2n+2+l(l-1)m^2+2lms+2m((l-1)m+s)=\]
\[ =-(ml+s)^2-2ml-2s+2+l(l-1)m^2+2lms+2m((l-1)m+s)=\]
\[ =(l-2)m^2+2ms-s^2-2(lm+s-1)=(l-2)m(m-2)+s(2m-s-2)-4m+2~.\]
One has $m\geq 2$, otherwise $\Lambda _2^n=(1,\ldots ,1)$. Hence, 
$s(2m-s-2)\geq 0$. 

$2^0$. If $l\geq 3$, $m\geq 6$ or $l\geq 4$, $m\geq 4$, 
then $(l-2)m(m-2)-4m>0$ and the triple cannot be rigid. On the other hand 
$l\geq 2$, otherwise $r_2\leq n/2$. Hence, rigid triples exist only for 
$l=2$ or 3 or for $m=2$ or 3. 

$3^0$. If $m=2$, then $\Delta <0$ (for $s=1$ or 2). The PMVs for which the 
minimal value of $\Delta$ is attained are:

\[ \begin{array}{llllllll}1) ~~(2,\ldots ,2)&,&(2,\ldots ,2)&,&(n-2,2)~&
{\rm for~}n~{\rm even}&;&\Delta =-6~;\\
2) ~~(2,\ldots ,2,1)&,&(2,\ldots ,2,1)&,&(n-2,2)~&{\rm for~}n~{\rm odd}&;
&\Delta =-4~.\end{array}\]
 
Find all rigid triples with such values of $r_2,r_3$ (i.e. $n-2,2$). To this 
end one has to replace 3 multiplicities equal to 2 for $n$ even 
(resp. 2 multiplicities equal to 2 for $n$ odd) by couples of multiplicities 
1,1. (Indeed, the biggest component of $\Lambda _1^n$ is $\leq 2$, the 
ones of $\Lambda _2^n$ 
and $\Lambda _3^n$ do not change because they define $r_2$ and $r_3$.) 
Each change of 2 by 1,1 increases $\Delta$ by 2. 

The possibilities (up to permutation of $\Lambda _1^n$ and 
$\Lambda _2^n$) for $n$ even are 1a) -- 1d), for $n$ odd they are 2a) -- 2c).
Possibility 2c) is the series $OG_k$ introduced in the previous section.

$4^0$. If $m=3$, $l\geq 5$, then $\Delta >0$, see $1^0$. Hence, for $m=3$  
rigid triples can exist only for $n\leq 15$. 

$5^0$. If $l=3$, then $\Delta =m(m-2)+s(2m-s-2)-4m+2$ and 
$\Delta >0$ if $m\geq 5$ or $m=4$, 
$s=2,3,4$ (to be checked directly). Hence, rigid triples with $l=3$ 
exist only for $n\leq 13$. 

$6^0$. If $l=2$, then $\Delta =s(2m-s-2)-4m+20$ and if $4\leq s\leq m-2$, 
then $\Delta >0$. 
Hence, rigid triples can exist only for $s=1,2,3,m-1,m$. 

If $s=3$ and $m\geq 7$, 
then $\Delta >0$, i.e. with $s=3$ rigid triples can exist only for $n\leq 21$. 

If $s=m$, then for $m\geq 6$ one has $\Delta >0$, i.e. such 
rigid triples can exist 
only for $n\leq 20$. 

If $s=m-1$, then again for $m\geq 6$ one has 
$\Delta >0$, i.e. such rigid triples can exist 
only for $n\leq 17$.  

If $s=1$, then we have $\Lambda _1^n=(2,\ldots ,2,1)$, $\Lambda _2^n=(m,m,1)$, 
$\Lambda _3^n=(m+1,m)$. One has $\Delta =-(n-1)$, i.e. for $n>1$ one cannot 
choose these MVs to have rigid triples. Give the list of the MVs with the same 
quantities $r_j$ for which $\Delta =0$. They are obtained from the given ones 
as a result of one or several antipassages, see Subsection~\ref{method}. 

The MV $\Lambda _2^n$ after one antipassage becomes $(m,m-1,2)$ (and $d_2$ 
increases by $(n-1)-4$) or $(m,m-1,1,1)$ (and $d_2$ increases by $(n-1)-2$). 
The MV $\Lambda _3^n$ after one antipassage becomes $(m+1,m-1,1)$ and $d_3$ 
increases by $(n-1)-2$. To increase $d_1$ by $2s$ one has to make $s$ 
antipassages 
in which a component 2 is replaced by a couple of units. However, we avoid 
to have $\Lambda _1^n=(1,\ldots ,1)$ which case was considered in 
Section~\ref{examples}. Therefore for $n\geq 22$ 
the only PMVs which give rigid triples for 
$s=1$ are $Z^2_n$, $Z^3_n$ and $Z^4_n$.  

If $s=2$, one gets the series $\Lambda _1^n=(2,\ldots ,2,2)$, 
$\Lambda _2^n=(m,m,2)$, $\Lambda _3^n=(m+2,m)$ with $\Delta =-6$. 
The only ways to increase $\Delta$ by 6 for $m\geq 10$ 
are to make three antipassages 
changing a component 2 by two components 1,1. This yields possibilities 1j) 
and 1k).

{\bf Subcase 1.2)} $r_2=r_3=n/2$ ($n$ is even).

We assume that $n\geq 22$. The PMV which minimizes the sum 
$d_1+d_2+d_3$ equals $(2,\ldots ,2)$, $(n/2,n/2)$, $(n/2,n/2)$ and one has 
$\Delta =-2n+2$. Hence, to obtain irreducible triples one has to choose 
another PMV, in which at least one MV defines a JNF giving a greater value 
of the corresponding quantity $d_j$. 

For the PMV as above one has $d_1=n^2-2n$, $d_2=d_3=n^2/2$. By replacing 
consecutively 
components equal to 2 of $\Lambda _1^n$ by couples of units one can obtain as 
values of $d_1$ all even numbers from $n^2-2n$ to $n^2-n$. 

Hence, one cannot increase enough 
$\Delta$ by changing only $\Lambda _1^n$. If one 
changes $\Lambda _2^n$ and/or $\Lambda _3^n$ without changing $r_2$ and $r_3$, 
the new choices have to be among the following MVs, otherwise $\Delta$ 
increases by more than $2n-2$:

\[ \begin{array}{lll}1)&(n/2,n/2-1,1)&d_j=n^2/2+n-2~;\\
2)&(n/2,n/2-2,2)&d_j=n^2/2+2n-8~;\\
3)&(n/2,n/2-2,1,1)&d_j=n^2/2+2n-6~.\end{array}\]
 
If one uses possibility 2) or 3), then the only cases in which $\Delta =0$ 
are 1g) and 1h). If one uses possibility 1), then this leads to case 1i) or to 
the series $EF_n$, see Section~\ref{examples}.

{\bf Case 2)} $r_2+r_3=n+1$. 

{\bf Subcase 2.1)} $r_2>n/2$, $r_3\leq n/2$.

The quantities $d_2$ and $d_3$ are minimal if and only if one 
has 

\[ \Lambda _2^n=(m,\ldots ,m,s)~,~\Lambda _3^n=(n-m-1,m+1)~,~
n=lm+s~,~l\in {\bf N}~,~1\leq s\leq m~,~m\geq 2~,~r_2=(l-1)m+s~,~r_3=m+1~,\] 
see Lemmas~\ref{min1} and \ref{min2}. 
For such $\Lambda _2^n$, $\Lambda _3^n$ 
one has 

\[ d_2=l(l-1)m^2+2lms~,~d_3=2(m+1)((l-1)m+s-1)=2m((l-1)m+s)+\delta \] 
where $\delta =2(l-1)m+2s-2-2m=2(l-2)m+2s-2$. Like in $1^0$ one finds 

\[ \Delta \geq (l-2)m(m-2)+s(2m-s-2)-4m+2+\delta =(l-2)m^2+s(2m-s)-4m\]
(the difference in the 
estimation of $d_3$ w.r.t. $1^0$ equals $\delta$). For $l>3$ one has 
$\Delta >0$. The 
same is true for $l=3$ except for $m=2$. In the latter case one has 
$n=7$ or 8.  

For $l=2$ one does not have $\Delta >0$ only if $s=1,2,3$ or 4; if 
$s=4$, then $m=s=4$ and $n=12$; if $s=3$, then $m=3$ or 4, resp. $n=9$ or 11. 

The case $l=2$, $s=1$ is impossible 
because then one has $r_2=m+1$, 
$r_3=m$ and $r_2+r_3=n<n+1$. 

If $l=2$, $s=2$, then $n=2m+2$ is even and for 
$\Lambda _1^n=(2,\ldots ,2)$, $\Lambda _2^n=(m,m,2)$, $\Lambda _3^n=(m+1,m+1)$ 
one has 

\[ \Delta =(n^2-2n)+2(n/2-1)^2+4(n-2)+n^2/2-2n^2+2=-4~.\]
If $m\leq 9$, then $n\leq 20$. If $m\geq 10$, i.e. $n\geq 22$, 
then it is possible to increase 
$\Delta$ by 4 (without changing $r_1$, $r_2$, $r_3$) only by replacing the PMV 
by one of the PMVs from 1e) or 1f).

For all other choices of $\Lambda _j^n$ with $r_1=n-2$, $r_2=m+2$, $r_3=m+1$ 
one has $\Delta >0$. Hence, no rigid triples exist for such PMVs.

{\bf Subcase 2.2)} $r_2>n/2$, $r_3>n/2$. 

Necessarily $n$ is odd ($n=2m+1$) and the 
minimal possible value of $d_1+d_2+d_3$ is attained for and only for  
$\Lambda _1^n=(2,\ldots ,2,1)$, $\Lambda _2^n=(m,m,1)$, $\Lambda _3^n=(m,m,1)$
(see Lemma~\ref{min2}). Such triples are rigid. They give possibility 2d).

Prove the existence of the listed series. With the notation from 
Section~\ref{examples} one has

\[ \Gamma ^i_n\rightarrow \Gamma ^i_{n-2}~~,~~\Gamma^1_6\rightarrow X_5^1
\rightarrow Y_4^1=\Gamma ^2_4=\Gamma ^4_4\rightarrow HG_3
~~,~~\Gamma^3_4\rightarrow HG_2~~,\]
which proves the existence of the series $\Gamma _n^i$, $i=1,2,3,4$. One 
also has 

\[ X_n^i\rightarrow X_{n-2}^i~~,~~X_5^1\rightarrow \Gamma ^2_4~~,~~
X^2_3=HG_3\] 
which proves the existence of the series $X^i_n$. Next, 

\[ {\rm for~~}n>4~~Y_n^1\rightarrow Z^2_{n-1}
\rightarrow Z^2_{n-3}\rightarrow \ldots \rightarrow Z^2_5
\rightarrow \Gamma _4^4~~,\]
hence, the series $Y_n^1$ and $Z^2_n$ also exist. From 
$Z_n^i\rightarrow Z_{n-2}^i$, $Z_3^i=HG_3$, $i=3,4$ there follows the 
existence of the series $Z_n^3$, $Z_n^4$. From 

\[ Z_n^1\rightarrow Y_{n-1}^5\rightarrow Y_{n-3}^5\rightarrow \ldots 
\rightarrow Y_2^5=HG_2\]
follows the existence of $Z_n^1$ and $Y_{n-1}^5$. From 
$Y_n^2\rightarrow Z_{n-1}^3$ follows the one of $Y_n^2$. One has 

\[ Y_n^3\rightarrow Y_{n-2}^6\rightarrow Y_{n-4}^3\rightarrow Y_{n-6}^6
\rightarrow \ldots ,\]
hence, $Y_n^6$ and $Y_n^3$ also exist (we let the reader prove the existence 
of $Y_4^3$; one has $Y_4^6=HG_4$). Finally, one has 

\[ Y_n^4\rightarrow Y_{n-2}^7\rightarrow Y_{n-4}^4\rightarrow Y_{n-6}^7
\rightarrow \ldots \]
which proves the existence of $Y_n^4$ and $Y_n^7$ (the reader has to prove the 
existence of $Y_6^4$ and $Y_6^7$).

The theorem is proved.

\section{The case of arbitrary (not necessarily diagonal) Jordan normal forms
\protect\label{nondiagonal}} 

\begin{defi}
{\em For a given JNF $J^n=\{ b_{i,l}\}$ define its {\em corresponding} 
diagonal JNF ${J'}^n$. (We say that $J^n$ and ${J'}^n$ are corresponding to 
one another.) A diagonal JNF is  
a partition of $n$ defined by the multiplicities of the eigenvalues. 
For each $l$ the family 
$\{ b_{i,l}\}$ is a partition of $\sum _{i\in I_l}b_{i,l}$ and 
${J'}^n$ is the disjoint sum of the dual partitions.}
\end{defi}

\begin{ex}
{\em If a JNF is defined by the family $B=\{ b_{i,l}\}$ where $l=1,2$ 
and $B=\{ 4,2,2\} \{5,1\}$, i.e. there are two eigenvalues, the first (resp. 
the second) with 
three Jordan blocks, of sizes 4, 2, 2 (resp. with two Jordan blocks, of sizes 
5, 1), then the corresponding diagonal JNF is defined by the MV 
$(3,3,1,1,4,2)$ (or, better, by the MV 
with non-increasing components $(4,3,3,2,1,1)$). Indeed, $(3,3,1,1)$ (resp. 
$(4,2)$) is the partition dual to $(4,2,2)$ (resp. to $(5,1)$).}
\end{ex}
 
The following theorem explains why it is sufficient to know (for generic 
eigenvalues) the solution to the DSP only in the case of diagonalizable 
matrices. The theorem is announced in \cite{Ko1} and proved in \cite{Ko2}. 

\begin{tm}\label{semisimple}
If for some eigenvalues the DSP is weakly solvable for a given 
$\{ J_j^n\}$ (resp. for $\{ {J'}_j^n\}$), then it is 
solvable for $\{ {J'}_j^n\}$ (resp. for $\{ J_j^n\}$) for 
any generic eigenvalues.
\end{tm}

Thus if one knows that the DSP is solvable for a certain PMV $\Lambda ^n$ 
for generic 
eigenvalues, then one knows that it is solvable (for generic eigenvalues) 
for all $(p+1)$-tuples of JNFs $\{ J_j^n\}$ such that the JNF defined by 
$\Lambda _j^n$ corresponds to $J_j^n$. This allows one to construct new series 
of $(p+1)$-tuples of JNFs (not all of which diagonal) 
for which there exist rigid $(p+1)$-tuples of matrices $A_j$ or $M_j$. 
One should know, however, that 
for certain  $(p+1)$-tuples of JNFs one cannot have generic eigenvalues.

\begin{ex}
{\em Consider the series $C_k$ from Section~\ref{examples} for matrices 
$A_j$. A possible triple 
of JNFs corresponding to the diagonal ones defined by the PMV 
is the following one: $J_1^n$ and $J_2^n$ are the same 
as before, i.e. diagonalizable, with MVs of the eigenvalues equal to 
$(k,k,k)$ while $J_j^3$ has a single eigenvalue with Jordan blocks of sizes 
$(1,2,3,\ldots ,3)$. Hence, the multiplicities of all eigenvalues are 
divisible by $k$. The sum of all eigenvalues counted with multiplicities $k$ 
times smaller equals 0 and this is a non-genericity relation.

Consider the same example for matrices $M_j$. The product of all eigenvalues 
with multiplicities $k$ times smaller is a root of unity of order $k$. If 
this root is 
non-primitive, then again a non-genericity relation holds and there exist 
no such generic eigenvalues. In this case the set of possible eigenvalues 
with these JNFs is a reducible variety with $k$ connected components 
each of which corresponds to one of the roots of unity. 
The eigenvalues from the components corresponding to non-primitive roots are 
all non-generic.}
\end{ex}

\section{Proof of Theorem~\ref{uleq2}\protect\label{pruleq2}}

$1^0$. Theorem \ref{semisimple} 
allows one to prove the theorem only in the case 
of diagonalizable matrices. For $n\leq 3$ the reader can check the theorem 
oneself, so suppose that $n\geq 4$. 

It suffices to prove that the PMV $\Lambda ^{n_1}$ obtained from 
$\Lambda ^n$ 
after applying $\Psi$ (see Section~\ref{DSPgen}) 
satisfies condition $(\beta _{n_1})$. (The PMV $\Lambda ^{n_1}$ satisfies 
condition $(\alpha _{n_1})$ if and only if 
$\Lambda ^n$ satisfies $(\alpha _n)$, see Lemma~\ref{invariant}.)

If one of the MVs is of the form 
$(1,\ldots ,1)$ and conditions $(\alpha _n)$, $(\beta _n)$ hold, then in the 
case of matrices $M_j$ the answer to the DSP is 
positive, see \cite{Si}, hence, it is positive for matrices $A_j$ as well (for 
generic eigenvalues the criterium is the same in the case of matrices $A_j$ 
or $M_j$). Therefore we assume that for all $j$ one has $m_{1,j}\geq 2$. 

\begin{rem}\label{dimension}
{\em Remind that 

1) the maximal value of $d_j$ equals $n^2-n$ and it is attained 
only for a MV of the form $(1,\ldots ,1)$;

2) for the MV $(n-1,1)$ the quantity $d_j$ equals $2n-2$; hence, if $p=2$ and 
one of the MVs equals $(n-1,1)$, then $(\alpha _n)$ holds only if the other 
two equal $(1,\ldots ,1)$;

3) for the MVs $(n/2,n/2)$ and $(n/2,n/2-1,1)$ the values of $d_j$ equal 
respectively $n^2/2$ and $n^2/2+n-2$.}
\end{rem}

$2^0$. Set $\rho _j:=r_1+\ldots +\hat{r}_j+\ldots +r_{p+1}$. 
One has $r_1=n-1$ or $n-2$, therefore for $j\neq 1$, $p\geq 3$ one has 
$\rho _j\geq n-2+p-1\geq n$; this is true for $p=2$ as well 
because no MV equals 
$(n-1,1)$, otherwise $(\alpha _n)$ does not hold, see Remark~\ref{dimension}. 
Therefore 
we check only that after performing the construction $\Psi$ from 
Section~\ref{DSPgen} one has $\rho _1\geq n_1$. 

One has $n-n_1\leq 2$, see Remark~\ref{diagcase}. If $n-n_1=1$, then every 
quantity $r_j$ remains the same or decreases by 1. The second possibility 
takes place only if $r_j\geq n/2$ and $\Lambda _j^n$ has two equal greatest 
components. Denote by $l$ the number of indices $j$ 
for which $r_j\geq n/2$. Hence, $j=1$ is always among them. 
Three cases are possible:

{\bf Case 1)} $l\leq 2$. 

Condition $(\beta _n)$ satisfied by $\Lambda ^n$ implies 
that $\Lambda ^{n_1}$ satisfies condition $(\beta _{n_1})$ because for $j>1$ 
either all $r_j$ remain the same or only one decreases by 1 when $\Psi$ 
is performed.  

{\bf Case 2)} $l\geq 3$ and $p\geq 3$. 

After applying $\Psi$ in the sum 
$\rho _1$ there are two quantities 
$r_j$ which are $\geq n/2-1$ and one which is $\geq 1$, so 
$\Lambda ^{n_1}$ satisfies condition 
$(\beta _{n_1})$. 

{\bf Case 3)}  $p=2$ and $l=3$.

The sum $\rho _1$  
can become $<n-1$ after applying $\Psi$ only if 
$\Lambda _2^n=\Lambda _3^n=(n/2,n/2)$ and 
$n$ is even. But in 
this case condition $(\alpha _n)$ does not hold for any MV $\Lambda _1^n$ 
(see Remark~\ref{dimension}), 
hence, the case has to be excluded. In all other cases the sum $\rho _1$ 
decreases by 1 and 
the PMV $\Lambda ^{n_1}$ satisfies condition 
$(\beta _{n_1})$. 

$3^0$. Let $n-n_1=2$. Like in the case $n-n_1=1$, for $j\neq 1$ the sum 
$\rho _j$ is $\geq n$. Indeed, if $p>2$, then such a sum contains 
$r_1\geq n-2$ and two more quantities $r_j$ which are $\geq 1$. If $p=2$ and 
$u=2$, then no MV is of the form $(n-1,1)$ because condition $(\alpha _n)$ 
would not hold, see Remark~\ref{dimension}. Hence, except 
$r_1\geq n-2$, $\rho _j$ contains $r_2\geq 2$ or 
$r_3\geq 2$, i.e. $\rho _j\geq n$. So there remains to check that after 
applying $\Psi$ one has $\rho _1\geq n-2$. 

$4^0$. Denote by $s_j$ the difference 
$m_{1,j}-m_{2,j}$ and 
by $\kappa$ the number of quantities $s_j$ which are $\leq 1$. Hence, $s_1$ 
is always one of them. Four cases are possible:

{\bf Case 4)} $\kappa =1$ or 2. 

At most one quantity $r_j$ from $\rho _1$ decreases by at most 2, so 
$\Lambda ^{n_1}$ satisfies condition $(\beta _{n_1})$. 

{\bf Case 5)} $\kappa \geq 4$. 

After performing $\Psi$ one has $r_j\geq n/2-2$ for three indices $j>1$, 
hence, $\rho _1\geq 3n/2-6\geq n-2$ because $n\geq 4$. 

{\bf Case 6)} $\kappa =3$, $p\geq 3$.

In this case after performing $\Psi$ one has $r_j\geq n/2-2$ for two 
indices $j>1$ and $r_j\geq 1$ for another one, so 
$\rho _1\geq n-3$ with 
equality only if two MVs $\Lambda _j^n$ with $j>1$ equal $(n/2,n/2)$ and 
a third equals $(n-1,1)$. But in such a case $n-n_1=1$, so the case has to 
be excluded. 

{\bf Case 7)} $\kappa =3$, $p=2$.

After performing $\Psi$ one has $\rho _1<n-2$ only if $n$ is even and 
either both 
$\Lambda _2^n$, $\Lambda _3^n$ are of the form $(n/2,n/2)$ or one is of this 
form while the other equals $(n/2,n/2-1,1)$. In the first case condition 
$(\alpha _n)$ does not hold for any $\Lambda _1^n$, see 
Remark~\ref{dimension}. In the second it holds 
only for $\Lambda _1^n=(1,\ldots ,1)$, but in this case $n-n_1=1$, so both 
cases have to be excluded.

The theorem is proved.

Author's address: Universit\'e de Nice, Laboratoire de Math\'ematiques, Parc 
Valrose, 06108 Nice Cedex 2, France; e-mail: kostov@math.unice.fr
\end{document}